\documentclass[11pt,a4paper]{article}
\usepackage[utf8]{inputenc}
\usepackage[T1]{fontenc}
\usepackage{amsmath,amssymb,amsthm}
\usepackage{mathpazo}
\usepackage[margin=2.5cm]{geometry}
\usepackage{mathrsfs}
\usepackage{booktabs}
\usepackage{microtype}
\usepackage{enumitem}
\usepackage{hyperref}
\usepackage{fancyhdr}
\usepackage{tikz}
\usetikzlibrary{arrows.meta, positioning}
\hypersetup{colorlinks=true,linkcolor=blue,urlcolor=blue,citecolor=blue}

\newtheorem{theorem}{Theorem}[section]
\newtheorem{proposition}[theorem]{Proposition}
\newtheorem{lemma}[theorem]{Lemma}

\newtheorem{definition}[theorem]{Definition}
\newtheorem{remark}[theorem]{Remark}
\newtheorem{example}[theorem]{Example}
\newtheorem{exercise}{Exercise}[section]

\newenvironment{keyidea}[1][Key idea:]
  {\medskip\noindent\begin{tabular}{|p{0.93\textwidth}|}
   \hline\\[-8pt]
   \textbf{#1} }
  {\\[4pt] \hline\end{tabular}\medskip}

\newcommand{\E}{\mathrm{E}}

\newcommand{\R}{\mathbb{R}}

\newcommand{\cS}{\mathcal{S}}
\newcommand{\cF}{\mathcal{F}}
\newcommand{\wm}[1]{{}^\varphi\!w_{#1}}
\newcommand{\wk}[1]{{}^\varphi\!\kappa_{#1}}

\newcommand{\wG}{{}^\varphi\!G}

\newcommand{\Ften} {\fontsize{10}{11}\selectfont  }

\title{Notes on Transversality and Statistical Degeneracies\\
in Distributional Models}

\author{Rodrigo Labouriau\\[4pt]
\small Department of Mathematics, Aarhus University\\
\small \texttt{rodrigo.labouriau@math.au.dk}\\
\small \texttt{rodrigo.labouriau@me.com}}
\date{(Spring 2026)}

\begin{document}
\maketitle
\thispagestyle{fancy}

\begin{abstract}
These notes provide a pedagogical introduction to the role of
transversality theory in the analysis of statistical degeneracies
within the framework of distributional statistical models.  The
classical question of when a statistical model is well-behaved---in
the sense of being identifiable, having non-singular Fisher
information, and admitting robust estimation---is reformulated as a
question about the geometry of a kernel-induced feature map.
Statistical pathologies correspond to geometric degeneracies of this
map, and transversality theory provides a precise language for
understanding when and why such degeneracies are non-generic.

The exposition is organised in three parts.  Part~I surveys the
statistical phenomena that motivate the geometric treatment:
representation failure, non-identifiability, moment indeterminacy,
singular information, nuisance parameters, and the Behrens--Fisher
problem.  Part~II develops the necessary geometric toolkit---smooth
maps, Sard's theorem, transversality, jets, stratifications, and the
parametric transversality theorem---at a level accessible to students
with a background in analysis and linear algebra but no prior
exposure to differential topology.  Part~III returns to the
statistical problems of Part~I and shows how each one admits a
unified geometric interpretation as a transversality condition on
the feature map.

These notes are a pedagogical companion to the research
paper~\cite{TransversalityPaper}, expanding its arguments with
motivating examples, geometric intuition, and exercises aimed at
advanced Master's and PhD students with a background in mathematical
statistics and measure theory.  They are designed to support seminars
or reading groups.
\end{abstract}

\newpage
\Ften
\tableofcontents
\normalsize
\newpage


\part*{Part I: Statistical Pathologies and the Distributional Framework}
\addcontentsline{toc}{part}{Part I: Statistical Pathologies and the
  Distributional Framework}

\bigskip

\section{Introduction and Motivation}
\label{sec:intro}

\textit{We introduce the central question of these notes: why do
  certain statistical models exhibit pathological behaviour, and what
  geometric structure underlies these pathologies?  We outline the
  distributional framework and explain how the introduction of a
  kernel resolves several difficulties simultaneously, raising the
  question that the rest of the notes will answer.}

\medskip

Statistical theory rests on a collection of regularity assumptions
that are so familiar as to be nearly invisible.  We assume that a
model is identifiable---that distinct parameter values give rise to
distinct distributions.  We assume that the Fisher information matrix
is non-singular, so that maximum likelihood estimators are
asymptotically efficient.  We assume that moments exist and determine
the distribution, so that method-of-moments estimators are
consistent.  We assume that densities exist with respect to a common
dominating measure, so that likelihood ratios are well-defined.

In many models of practical and theoretical importance, one or more of
these assumptions fail.  The log-normal distribution is
moment-indeterminate: infinitely many distinct distributions share the
same sequence of moments.  The Cauchy distribution has no finite
moments at all: the classical moment map is undefined.  Elliptically
contoured distributions defined through their characteristic function
may lack closed-form densities.  In the Behrens--Fisher problem, no
exact test exists because the null hypothesis is not in a ``generic
position'' relative to the nuisance parameters.

The purpose of these notes is to develop a geometric perspective on
these phenomena.  We will see that each of the pathologies listed above
corresponds to a specific geometric degeneracy of a map associated
with the statistical model, and that the theory of transversality
provides a unified framework for understanding when such degeneracies
occur and how they can be resolved.

\begin{keyidea}
Statistical pathologies---non-identifiability, singular information,
moment indeterminacy, representation failure---are geometric
degeneracies of a feature map associated with the model.
Transversality theory explains why they are non-generic and how
a kernel resolves them.
\end{keyidea}

The geometric perspective is made possible by a recent line of work
that replaces classical probability densities by
\emph{distribution--kernel pairs} $(T,\varphi)$, where $T$ is a
tempered distribution and $\varphi$ is a rapidly decaying kernel
\cite{A}.  In this framework, expectations are defined through
distributional pairings rather than through integration against a
density, and moments are replaced by \emph{weak moments}:
\[
  \wm{j}(\theta)
  \;=\;
  \langle T_\theta,\, x^j\,\varphi(x)\rangle.
\]
These weak moments are well-defined for all orders and all models,
because the rapid decay of the kernel provides the necessary
integrability.  The kernel thus induces a \emph{feature map}
\[
  \Phi_\varphi : \Theta \to \cF,
  \qquad
  \Phi_\varphi(\theta)
  = \bigl(\wm{0}(\theta), \wm{1}(\theta), \wm{2}(\theta),
  \ldots\bigr),
\]
which maps the parameter space into a space of weak moments.  The
statistical properties of the model are encoded in the geometry of
this map.

Empirically, the introduction of a kernel resolves the pathologies
listed above: the log-normal becomes moment-determinate, the Cauchy
acquires finite weak moments and regular information, influence
functions become bounded, and the Behrens--Fisher problem admits an
approximate solution.  The question left open---and the question
these notes aim to answer---is:

\begin{keyidea}[Central question:]
Why does the introduction of a kernel generically resolve statistical
degeneracies?
\end{keyidea}

The answer, as we shall see, is that the kernel acts as a
\emph{generic perturbation} in the sense of transversality theory.  The
parametric transversality theorem guarantees that for ``most'' kernels
in a sufficiently rich family, the feature map avoids degeneracy
strata, just as a generic smooth map avoids the singular set of its
target.

The exposition is organised as follows.  Sections~\ref{sec:representation}
through~\ref{sec:behrens_fisher_stat} survey the statistical phenomena
that motivate the geometric treatment, and
Section~\ref{sec:distributional_framework} introduces the distributional
framework.  Part~II (Sections~\ref{sec:smooth_maps}
through~\ref{sec:infinite_dim}) develops the geometric toolkit.
Part~III (Sections~\ref{sec:degeneracy_strat}
through~\ref{sec:singular_limit}) returns to the statistical problems
and provides their geometric interpretation.

\section{Representing Probability Distributions}
\label{sec:representation}

\textit{We discuss the classical representation of probability
  distributions via Radon--Nikodym densities, its dependence on the
  choice of dominating measure, and models that lack closed-form
  densities altogether.  This motivates the introduction of
  distribution--kernel pairs.}

\medskip

The standard approach to parametric statistics begins with a family of
probability measures $\{P_\theta : \theta \in \Theta\}$ on a sample
space $(\mathcal{X}, \mathcal{A})$, where $\Theta \subset \R^p$ is an
open set.  If the family is dominated by a $\sigma$-finite measure
$\mu$, the Radon--Nikodym theorem provides densities
$f(x;\theta) = dP_\theta/d\mu(x)$, and the entire machinery of
likelihood-based inference can be deployed.

However, this representation involves a choice: the dominating measure
$\mu$ is not unique, and the density $f(x;\theta)$ is only determined
up to $\mu$-null sets.  In most textbook treatments, this ambiguity is
harmless.  But it becomes problematic in at least two ways.

First, the choice of $\mu$ can affect the apparent regularity of the
model.  Different representatives of the density may have different
smoothness properties, and quantities such as the score function
$\partial_\theta \log f(x;\theta)$ depend on the choice of
representative at each point~$x$.

Second, some models of interest do not admit closed-form densities at
all.  Consider a family of elliptically contoured distributions
defined by the characteristic function
$\varphi_X(u) = e^{iu^\top\mu}\psi(u^\top\Sigma u)$, where $\psi$
does not admit closed-form Fourier inversion.  The distributions are
perfectly well-defined as probability measures, but the classical
density-based framework cannot even be set up.

\begin{keyidea}
The classical representation of a statistical model via densities
depends on an arbitrary choice (the dominating measure) and may not
exist at all.  This is a \emph{representation degeneracy}: the model
cannot be embedded into the standard coordinate system.
\end{keyidea}

\begin{remark}[Statistical perspective]
\label{rem:representation_stat}
From a practical viewpoint, representation degeneracy prevents the use
of likelihood-based methods: if no density is available, the
likelihood function cannot be written down.  This affects not only
estimation but also model comparison, hypothesis testing, and
information-theoretic quantities.
\end{remark}

\begin{remark}[Geometric perspective]
\label{rem:representation_geom}
From a geometric viewpoint, representation degeneracy means that the
feature map $\Phi_\varphi$ is \emph{undefined} in the classical sense.
In the distributional framework, the kernel creates an embedding: the
pairing $\langle T_\theta, g_j\varphi\rangle$ is well-defined and
smooth even when no density exists.  The kernel thus resolves the
representation degeneracy by providing a coordinate system that does
not depend on the existence of a density
(see~\cite{TransversalityPaper}, Section~7, Type~0 degeneracy and
Example~8.3 for the elliptically contoured case).
\end{remark}

\section{Identifiability and the Moment Problem}
\label{sec:identifiability}

\textit{We discuss identifiability as the requirement that the
  parametrisation of a model be injective, and connect it to the
  classical moment problem.  The log-normal distribution provides the
  central example of a model that is identifiable in the parametric
  sense but moment-indeterminate in the distributional sense.}

\medskip

A parametric model $\{P_\theta : \theta \in \Theta\}$ is
\emph{identifiable} if the map $\theta \mapsto P_\theta$ is injective:
distinct parameter values give rise to distinct distributions.
Identifiability is the minimal requirement for consistent estimation.
If the model is not identifiable, the data cannot distinguish between
different parameter values, and the statistical problem is
fundamentally ill-posed.

A closely related question, arising in the method of moments, is
whether the moment sequence determines the distribution.  Given a
probability measure $\mu$ on $\R$ with moments of all orders,
$m_n = \int x^n \, d\mu(x)$, the \emph{moment problem} asks: does
the sequence $(m_n)_{n \geq 0}$ uniquely determine $\mu$?

A distribution is said to be \emph{M-determinate} if the answer is
yes, and \emph{M-indeterminate} otherwise.  When the moment problem
has a positive answer, identifiability from moments is guaranteed.

\paragraph{A motivating example.}
The standard log-normal distribution provides a striking example of
M-indeterminacy.  Let $\mu$ denote the law of $e^Y$ with
$Y \sim N(0,1)$.  All moments exist and are given by
$m_n = \exp(n^2/2)$.  A classical construction due to Heyde produces a
one-parameter family of probability measures $\mu_a$,
$a \in [-1,1]$, all with the \emph{same} moment sequence as $\mu$:
\[
  \mu_a(dx) = \bigl(1 + a\sin(2\pi\log x)\bigr)\,\mu(dx),
  \qquad x > 0.
\]
Each $\mu_a$ is a probability measure, and
$\int x^n \sin(2\pi \log x)\,\mu(dx) = 0$ for every $n$ (by a direct
calculation using the symmetry of the Gaussian).  Thus the moment
sequence $(m_n)$ does \emph{not} determine the distribution.

The underlying reason is that the log-normal moments grow too rapidly:
$m_n = \exp(n^2/2)$ grows faster than any exponential in $n$.  This
violates the Carleman condition
$\sum_{n=1}^\infty m_{2n}^{-1/(2n)} = \infty$, which is a sufficient
condition for M-determinacy (see \cite{QuasiAnalNotes} for a detailed
treatment via the Denjoy--Carleman framework).

\begin{keyidea}
M-indeterminacy means that the moment map fails to be injective on
the space of distributions.  In geometric language, the feature map
has self-intersections: distinct distributions are mapped to the same
point in feature space.
\end{keyidea}

\begin{remark}[Statistical perspective]
\label{rem:M-indet_stat}
From a statistical viewpoint, M-indeterminacy undermines moment-based
estimation.  If the moment map is not injective, a method-of-moments
estimator may converge to a parameter value that does not correspond
to the true distribution.  The Stieltjes perturbation
$h(x) = \sin(2\pi \log x)$ represents a direction in distribution
space along which the moment map is completely blind.
\end{remark}

\begin{remark}[Geometric perspective]
\label{rem:M-indet_geom}
From a geometric viewpoint, M-indeterminacy is a failure of
transversality of the moment map relative to the self-intersection
diagonal.  The Stieltjes perturbation $h(x) = \sin(2\pi \log x)$ is
the tangent direction along which the map degenerates.  As we will see
in Part~III, a Gaussian kernel breaks the theta-function symmetry
underlying this cancellation, restoring transversality and hence
injectivity for generic kernels
(see~\cite{TransversalityPaper}, Section~8.1 and
Proposition~8.1 for the full analysis).
\end{remark}

\section{Fisher Information and Its Singularities}
\label{sec:fisher_info}

\textit{We introduce the Fisher information matrix as a measure of the
  information content of a statistical model, discuss its role in
  asymptotic theory, and describe situations in which it degenerates.
  This motivates the geometric interpretation of information regularity
  as an immersion condition on the feature map.}

\medskip

Given a parametric model with density $f(x;\theta)$, the \emph{score
  function} is $s(x;\theta) = \nabla_\theta \log f(x;\theta)$, and the
\emph{Fisher information matrix} is
\[
  I(\theta)_{ab}
  = \E_\theta\!\left[
    \frac{\partial \log f}{\partial \theta_a}\,
    \frac{\partial \log f}{\partial \theta_b}
  \right].
\]
Under standard regularity conditions, the Fisher information governs
the asymptotic variance of the maximum likelihood estimator via the
Cram\'er--Rao bound: no unbiased estimator can have variance smaller
than $I(\theta)^{-1}/n$.

When the Fisher information matrix is non-singular, the model is said
to be \emph{information-regular}.  But singularity of $I(\theta)$ can
occur in several ways.

\paragraph{Heavy tails.}
The Cauchy location family
$f(x;\mu) = \frac{1}{\pi}\frac{1}{1+(x-\mu)^2}$ has no finite
moments.  The score function $s(x;\mu) = 2(x-\mu)/(1+(x-\mu)^2)$
does have finite variance, so the Fisher information exists and equals
$I(\mu) = 1/2$.  But this is a special case; for many heavy-tailed
distributions, the score may not be square-integrable, and the Fisher
information may be infinite or undefined.

\paragraph{Models without densities.}
If no density exists (Section~\ref{sec:representation}), the score
function is undefined, and the Fisher information cannot even be
formulated.

\paragraph{Boundary effects and curved models.}
At the boundary of the parameter space, or in models with constraints,
the Fisher information may drop rank.  This occurs, for instance, in
mixture models at the boundary where a component weight vanishes.

\begin{keyidea}
The Fisher information matrix measures the local geometry of the
model: it is the metric tensor of the statistical manifold.
Singularity of the information matrix corresponds to the feature map
failing to be an immersion---the model ``collapses'' in certain
parameter directions.
\end{keyidea}

In the distributional framework, the classical Fisher information is
replaced by the \emph{distributional information tensor}
\[
  \wG_{ab}^{(J)}(\theta)
  = \sum_{j=0}^{J}
  \frac{\partial \wm{j}}{\partial\theta_a}\,
  \frac{\partial \wm{j}}{\partial\theta_b},
\]
which is the first fundamental form of the feature-map immersion
$\Phi_\varphi$.  This tensor is well-defined and smooth for any model,
including those without densities, because the kernel provides the
necessary regularity.

\begin{remark}[Statistical perspective]
\label{rem:info_stat}
From a statistical viewpoint, singularity of the information matrix
means that the model is locally flat in some direction: infinitesimal
changes in the parameter do not produce detectable changes in the
distribution.  This makes consistent estimation impossible along the
degenerate direction.
\end{remark}

\begin{remark}[Geometric perspective]
\label{rem:info_geom}
From a geometric viewpoint, information singularity is a
Thom--Boardman singularity of the feature map: the derivative
$D\Phi_\varphi(\theta)$ drops rank.  Transversality theory predicts
that such rank-drop loci are generically of positive codimension,
and hence avoidable by a suitable perturbation.  The kernel provides
exactly this perturbation
(see~\cite{TransversalityPaper}, Section~7, Type~II degeneracy,
and Corollary~4.2 for the generic full-rank result).
\end{remark}

\section{Robustness and Influence Functions}
\label{sec:robustness}

\textit{We discuss the sensitivity of statistical procedures to
  outliers and model misspecification, introduce the influence function
  as a measure of robustness, and explain why classical procedures for
  heavy-tailed models are inherently non-robust.}

\medskip

A statistical procedure is \emph{robust} if its performance degrades
gracefully under small departures from the assumed model.  The
classical tool for measuring robustness is the \emph{influence
  function}: given a functional $T(F)$ of the distribution function
$F$, the influence function at a point $x$ is
\[
  \mathrm{IF}(x; T, F)
  = \lim_{\varepsilon \to 0}
  \frac{T\bigl((1-\varepsilon)F + \varepsilon\delta_x\bigr) - T(F)}{\varepsilon},
\]
where $\delta_x$ is the point mass at~$x$.  A bounded influence
function indicates that no single observation can have an arbitrarily
large effect on the estimator.

For heavy-tailed distributions, classical estimators (such as the
sample mean) have unbounded influence functions: a single outlier can
move the estimator arbitrarily far.  This is a practical manifestation
of the failure of moment conditions.

In the distributional framework, the kernel provides a natural
damping mechanism.  The weak moment
$\wm{j}(\theta) = \E_\theta[X^j \varphi(X)]$ involves the factor
$\varphi(X)$, which decays rapidly for large $|X|$.  This ensures
that the influence of extreme observations is bounded by the decay
rate of the kernel, leading to estimators with bounded influence
functions even for heavy-tailed models.

\begin{keyidea}
Robustness corresponds to boundedness of the distributional metric:
the feature-map immersion has bounded first fundamental form.  The
kernel ensures this by dampening the contribution of extreme
observations.
\end{keyidea}

\section{Nuisance Parameters and the Behrens--Fisher Problem}
\label{sec:behrens_fisher_stat}

\textit{We introduce the problem of inference in the presence of
  nuisance parameters, describe the Behrens--Fisher problem as a
  classical instance where exact inference is impossible, and connect
  inferential separation to the geometry of projections along
  nuisance directions.}

\medskip

In many statistical problems, the parameter $\theta$ decomposes as
$\theta = (\psi, \phi)$, where $\psi$ is the \emph{parameter of
  interest} and $\phi$ is a \emph{nuisance parameter}.  The goal is to
make inference about $\psi$ without being affected by the unknown
value of $\phi$.

When exact separation of interest and nuisance is possible---for
instance, through a sufficient statistic for $\phi$ that is
ancillary for $\psi$---the problem reduces to inference in a
lower-dimensional model.  This ideal situation is called
\emph{inferential separation}.

The Behrens--Fisher problem is a classical instance where inferential
separation fails.  Let $X_1, \ldots, X_m \sim N(\mu_1, \sigma_1^2)$
and $Y_1, \ldots, Y_n \sim N(\mu_2, \sigma_2^2)$ be independent, with
$\sigma_1^2 \neq \sigma_2^2$ both unknown.  The problem is to test
$H_0: \mu_1 = \mu_2$.  No pivotal quantity exists whose distribution
is free of the nuisance ratio $\rho = \sigma_1^2/\sigma_2^2$: the
distribution of any test statistic depends on $\rho$, and exact
inference about the means is impossible without knowing the variance
ratio.

\begin{keyidea}
The Behrens--Fisher problem is a failure of the classical feature map
to place the null hypothesis in a generic position relative to the
nuisance parameter structure.  The distribution of any test statistic
depends on the nuisance parameters because the testing direction is
not transversal to the nuisance fibration.
\end{keyidea}

More broadly, inferential separation can be understood through the
\emph{nonformation principle}, which provides a hierarchy of
increasingly weak separation conditions:
\emph{S-nonformation} (the feature map factors completely through the
interest projection), \emph{I-nonformation} (local separation at the
maximum), and \emph{L-nonformation} (separation at the level of the
normed profile).  As we will see in Part~III, these correspond to
progressively weaker transversality conditions on the feature map
relative to the nuisance-parameter fibration.

\section{The Distributional Framework}
\label{sec:distributional_framework}

\textit{We introduce the distributional statistical framework:
  distribution--kernel pairs, weak moments, weak characteristic
  functions, and weak cumulants.  We explain how the kernel resolves
  each of the pathologies discussed in the preceding sections, and
  introduce the feature map that will be the central object of the
  geometric analysis.}

\medskip

The distributional statistical framework, developed in \cite{A,B}
and summarised in~\cite{TransversalityPaper}, Section~2,
replaces the classical density by a pair $(T_\theta, \varphi)$, where
$T_\theta \in \cS'(\R^d)$ is a tempered distribution (in the sense of
Schwartz) parametrised by $\theta \in \Theta \subset \R^p$, and
$\varphi \in \cS(\R^d)$ is a positive Schwartz function called the
\emph{kernel}.  Expectations are defined through distributional
pairings:
\[
  {}^{(\varphi)}\E_\theta[g]
  = \langle T_\theta, g\varphi \rangle.
\]

\paragraph{Weak moments.}
The \emph{weak moment of order~$j$} is
\[
  \wm{j}(\theta)
  = \langle T_\theta, x^j \varphi(x) \rangle
  = \E_\theta[X^j \varphi(X)].
\]
Because $\varphi$ is rapidly decaying, the product $x^j \varphi(x)$ is
integrable for every $j$, and weak moments of all orders are
well-defined for any distribution---including the Cauchy, stable laws,
and distributions without densities.

\paragraph{Weak characteristic function.}
The \emph{weak characteristic function} is
${}^{(\varphi)}\phi_\theta(u) = \langle T_\theta,
e^{iux}\varphi(x)\rangle$.  It is an entire function of~$u$, even
when the classical characteristic function has limited smoothness.

\paragraph{Weak cumulants.}
The kernel defines a tilted probability
$p_\varphi(x;\theta) = f(x;\theta)\varphi(x)/\wm{0}(\theta)$.
The \emph{weak cumulants} $\wk{j}$ are the cumulants of $X$ under
this tilted distribution.  The weak cumulant generating function
is entire, even when the classical cumulant generating function does
not exist.

\paragraph{The feature map.}
Given a kernel $\varphi$ and a collection of moment orders
$0 \leq j_0 < j_1 < \cdots < j_K$, the \emph{feature map} is
\[
  \Phi_\varphi : \Theta \to \R^{K+1},
  \qquad
  \Phi_\varphi(\theta)
  = \bigl(\wm{j_0}(\theta), \ldots, \wm{j_K}(\theta)\bigr).
\]
This map encodes the statistical properties of the model:
identifiability corresponds to injectivity, information regularity to
immersivity, and robustness to boundedness of the induced metric.

\paragraph{How the kernel resolves pathologies.}
The pathologies discussed in the preceding sections are resolved as
follows:

\begin{center}
\begin{tabular}{lll}
\toprule
\textbf{Pathology} & \textbf{Classical failure} &
  \textbf{Kernel resolution} \\
\midrule
Representation & No density exists & Pairing $\langle T_\theta,
  g\varphi\rangle$ is well-defined \\
Identifiability & Moment map not injective & Weak moment map
  generically injective \\
Information & $I(\theta)$ singular or undefined &
  $\wG(\theta)$ well-defined and generic.\ non-singular \\
Robustness & Unbounded influence & Kernel decay bounds influence \\
Behrens--Fisher & Nuisance non-separation & Kernel couples
  $\mu$ and $\sigma$ via $\sigma^2+s^2$ \\
\bottomrule
\end{tabular}
\end{center}

\medskip

The question posed in the introduction can now be made precise:

\begin{keyidea}[Precise formulation of the central question:]
Under what conditions does the feature map $\Phi_\varphi$ avoid the
degeneracy loci (self-intersections, rank drops, etc.) for a generic
choice of kernel $\varphi$?
\end{keyidea}

The answer will be provided by the parametric transversality theorem,
developed in Part~II.

\paragraph{A remark on the scope of the present notes.}
The transversality results developed here apply to \emph{parametric}
models, i.e.\ models with finite-dimensional parameter spaces
$\Theta \subset \R^p$.  The extension to semiparametric and
nonparametric settings, where $\Theta$ is infinite-dimensional,
requires the Banach-manifold framework of Smale and Quinn.  We discuss
the conceptual picture for these extensions in
Section~\ref{sec:infinite_dim}, but the full development is deferred
to future work.


\newpage
\part*{Part II: Transversality --- A Geometric Toolkit}
\addcontentsline{toc}{part}{Part II: Transversality --- A Geometric Toolkit}

\bigskip

\section{Smooth Maps and Their Singularities}
\label{sec:smooth_maps}

\textit{We introduce smooth maps between Euclidean spaces, the
  Jacobian, regular and critical values, and state Sard's theorem.  The
  preimage theorem shows that regular level sets are smooth manifolds.
  These ideas provide the foundation for the transversality theory
  developed in the following sections.}

\medskip

Let $f : \R^m \to \R^n$ be a smooth ($C^\infty$) map.  At each point
$x \in \R^m$, the \emph{Jacobian} (or derivative) of $f$ is the
linear map
\[
  Df(x) : \R^m \to \R^n,
\]
represented by the $n \times m$ matrix of partial derivatives
$(\partial f_i / \partial x_j)$.  The \emph{rank} of $f$ at $x$ is
the rank of this matrix.

\begin{definition}[Regular and critical values]
\label{def:regular_value}
A point $y \in \R^n$ is a \emph{regular value} of $f$ if for every
$x \in f^{-1}(y)$, the Jacobian $Df(x)$ is surjective (i.e.\ has
rank $n$).  A point that is not a regular value is a \emph{critical
  value}.  The set of critical values is the \emph{critical set} of
$f$.
\end{definition}

The fundamental result about critical values is Sard's theorem, which
asserts that critical values are ``negligible'' in a precise sense.

\begin{theorem}[Sard, 1942]
\label{thm:sard}
Let $f : \R^m \to \R^n$ be a smooth map.  Then the set of critical
values of $f$ has Lebesgue measure zero in $\R^n$.
\end{theorem}

In other words, ``almost every'' value of $f$ is regular.  This is a
remarkable result: it says that the generic behaviour of a smooth map
is to be surjective at every preimage, regardless of how complicated
the map may be.

The importance of regular values lies in the following consequence:

\begin{theorem}[Preimage theorem]
\label{thm:preimage}
If $y$ is a regular value of $f : \R^m \to \R^n$ with $m \geq n$,
then $f^{-1}(y)$ is either empty or a smooth submanifold of $\R^m$ of
dimension $m - n$.
\end{theorem}

Geometrically, this means that the level sets of a smooth map are
``well-behaved'' (smooth submanifolds) at generic values, and can only
develop singularities at the critical values---which form a negligible
set.

\begin{example}
\label{ex:sphere}
The map $f : \R^3 \to \R$ defined by $f(x,y,z) = x^2 + y^2 + z^2$
has $Df(x,y,z) = (2x, 2y, 2z)$, which is surjective whenever
$(x,y,z) \neq (0,0,0)$.  Thus every $r > 0$ is a regular value, and
$f^{-1}(r)$ is a smooth $2$-sphere.  The only critical value is
$r = 0$, where the level set degenerates to a point.
\end{example}

\begin{keyidea}
Sard's theorem says that the singular behaviour of a smooth map is
non-generic: for almost every target value, the preimage is a smooth
manifold of the expected dimension.  This is the prototype for the
transversality results that follow.
\end{keyidea}

\begin{remark}[Statistical perspective]
\label{rem:sard_stat}
In the statistical context, the map $f$ will be the feature map
$\Phi_\varphi$, the domain will be the parameter space $\Theta$, and
the target will be the feature space of weak moments.  Sard's theorem
tells us that, generically, the feature map behaves well: its level
sets (the sets of parameters giving rise to the same weak moments)
are smooth submanifolds---or, ideally, single points (when the map is
injective).
\end{remark}

\begin{remark}[Geometric perspective]
\label{rem:sard_geom}
Sard's theorem is the starting point of singularity theory: it tells
us that generic maps have controlled singularities.  Transversality
theory (next section) refines this by specifying exactly how a map
can be ``perturbed'' into a generic position.  For a concise treatment
in the statistical context, see~\cite{TransversalityPaper}, Section~3.
\end{remark}

\section{Transversality}
\label{sec:transversality}

\textit{We define transversality of a smooth map to a submanifold,
  explain its geometric meaning through examples, and state the
  fundamental consequence: transversal preimages are smooth
  submanifolds of the expected codimension.}

\medskip

Sard's theorem deals with the special case where the target
submanifold is a single point.  Transversality generalises this to
arbitrary submanifolds.

Let $M$ and $N$ be smooth manifolds (for concreteness, open subsets of
Euclidean spaces), and let $S \subset N$ be a smooth submanifold of
codimension $c$.

\begin{definition}[Transversality]
\label{def:transversality}
A smooth map $f : M \to N$ is \emph{transversal} to $S$, written
$f \pitchfork S$, if for every $x \in f^{-1}(S)$,
\[
  Df_x(T_x M) + T_{f(x)} S = T_{f(x)} N.
\]
In words: the image of the derivative, together with the tangent
space of $S$, spans the full tangent space of $N$.
\end{definition}

The condition says that $f$ meets $S$ in a ``non-degenerate'' way: the
map approaches $S$ from sufficiently many independent directions to
fill out the ambient space.

\paragraph{Geometric illustrations.}

Consider a line and a surface in $\R^3$.  The line is transversal to
the surface if it crosses the surface at an angle (not tangentially).
If the line is tangent to the surface at the intersection point, the
transversality condition fails: the image of the derivative of the
line (a one-dimensional subspace) lies entirely in the tangent plane
of the surface, and together they do not span $\R^3$.

More formally, consider a curve
$\gamma : \R \to \R^2$ and a curve $S \subset \R^2$.  Transversality
$\gamma \pitchfork S$ at a point $\gamma(t_0) \in S$ means that the
velocity vector $\gamma'(t_0)$ is not tangent to $S$: the two curves
cross rather than touch.

\paragraph{The preimage theorem for transversal maps.}

\begin{proposition}
\label{prop:transversal_preimage}
If $f \pitchfork S$, then $f^{-1}(S)$ is either empty or a smooth
submanifold of $M$ of codimension $c$ (the same codimension as $S$ in
$N$).
\end{proposition}

This is the transversal generalisation of Theorem~\ref{thm:preimage}.
In particular, when $c > \dim M$, transversality forces
$f^{-1}(S) = \emptyset$: the image of $f$ misses $S$ entirely.

\begin{keyidea}
Transversality is the ``right'' notion of genericity for the
intersection of a map with a submanifold.  Transversal intersections
have the expected codimension; non-transversal intersections are
``accidental'' and can be removed by a small perturbation.
\end{keyidea}

\begin{keyidea}[Codimension counting:]
If $f : M \to N$ is transversal to $S \subset N$ and
$\mathrm{codim}(S) > \dim M$, then $f(M) \cap S = \emptyset$:
the map avoids $S$ entirely.  This is the mechanism by which
sufficiently many weak moments guarantee identifiability.
\end{keyidea}

\section{Genericity: Thom's Transversality Theorem}
\label{sec:thom}

\textit{We state Thom's transversality theorem, which asserts that
  transversality is a generic condition, and the parametric version,
  which shows that a family of maps can be ``tuned'' into a
  transversal position by varying a parameter.  The parametric version
  is the key result for the statistical application: the kernel plays
  the role of the parameter.}

\medskip

The examples of the previous section might suggest that transversality
is a special condition, difficult to achieve.  In fact, the opposite
is true: transversality is the \emph{generic} condition, and
non-transversality is the exception.

\begin{theorem}[Thom's Transversality Theorem]
\label{thm:thom}
Let $M$ and $N$ be smooth manifolds and $S \subset N$ a smooth
submanifold.  Then the set
\[
  \{f \in C^\infty(M,N) : f \pitchfork S\}
\]
is residual (a countable intersection of open dense sets) in the
Whitney $C^\infty$ topology.  If $S$ is closed, this set is also
open dense.
\end{theorem}

In plain language: ``almost every'' smooth map is transversal to any
given submanifold.  Non-transversal maps are exceptional, like the
rational numbers among the reals---they are everywhere, but they form
a negligible set.

For the statistical application, we need a more structured result.  We
do not want to perturb the feature map arbitrarily; we want to perturb
it by varying the kernel $\varphi$, which enters as a parameter.  This
leads to the parametric version.

\begin{theorem}[Parametric Transversality Theorem]
\label{thm:parametric}
Let $F : M \times \Lambda \to N$ be a smooth map, where $\Lambda$ is
a smooth manifold of parameters.  Suppose $F \pitchfork S$.  Then
for a residual set of $\lambda \in \Lambda$, the restricted map
$F_\lambda = F(\cdot, \lambda) : M \to N$ satisfies
$F_\lambda \pitchfork S$.
\end{theorem}

The logic is as follows.  The full map $F$ depends on both the
``intrinsic'' variable $x \in M$ and the parameter
$\lambda \in \Lambda$.  If the full map is transversal to $S$
(which is easier to verify, since $F$ has more degrees of
freedom), then for \emph{generic} values of the parameter, the
restricted map is also transversal.

\begin{keyidea}
The parametric transversality theorem is the engine of the entire
theory.  In the statistical application, $M = \Theta$ (parameter
space), $\Lambda$ is the kernel family, $N$ is the feature space, and
$S = D$ is a degeneracy stratum.  The theorem says: if the joint
feature map $F(\theta, \lambda) = \Phi_\lambda(\theta)$ is
transversal to $D$, then for a generic kernel, the feature map avoids
$D$ (or meets it in the expected codimension).
\end{keyidea}

\begin{remark}[Statistical perspective]
\label{rem:parametric_stat}
The kernel plays the role of the parameter $\lambda$.  The result says
that if the model and kernel family are ``jointly rich enough'' (the
full map $F$ is transversal), then for most choices of kernel, the
feature map is well-behaved.  The non-transversal kernels---including
the degenerate limit $\varphi \to 1$ (classical moments)---form a
negligible set.
\end{remark}

\begin{remark}[Geometric perspective]
\label{rem:parametric_geom}
The parametric transversality theorem reduces the problem of
establishing genericity to a single ``global'' transversality
verification.  In practice, this is often done by checking that the
Jacobian of the joint map $F$ has full rank---the verifiable conditions
developed in~\cite{TransversalityPaper}, Section~5 and revisited
in Part~III of these notes.
\end{remark}

\section{Jets, Multijets, and Stratifications}
\label{sec:jets}

\textit{We introduce jets as the formal Taylor polynomials of smooth
  maps, explain how different types of singularity correspond to
  conditions on jets of different orders, and describe Whitney
  stratifications as a tool for organising singular sets.  The
  Thom--Mather theory classifies the generic singularities of smooth
  maps.}

\medskip

The transversality theorem of the previous section guarantees that a
generic map avoids a given submanifold (or meets it transversally).
But in practice, the degeneracy loci of a feature map are not single
submanifolds: they form a hierarchy of strata, corresponding to
different types and severities of degeneracy.  To describe this
hierarchy, we need the language of jets.

\paragraph{Jets.}
The \emph{$r$-jet} of a smooth map $f : M \to N$ at a point $x$ is
the equivalence class of $f$ modulo maps that agree with $f$ to order
$r$ at $x$.  Concretely, it is the collection of Taylor coefficients
up to order $r$:
\[
  j^r f(x) = \bigl(x, f(x), Df(x), D^2f(x), \ldots, D^rf(x)\bigr).
\]
The space of all $r$-jets forms a smooth manifold $J^r(M,N)$, called
the \emph{$r$-jet space}.  The map $j^r f : M \to J^r(M,N)$ that
sends each point to its $r$-jet is called the \emph{$r$-jet
  extension} of $f$.

\paragraph{Why jets matter.}
Different types of singularity correspond to conditions on jets of
different orders:

The $0$-jet $j^0 f(x) = (x, f(x))$ records the value of $f$.  Two
points $x_1 \neq x_2$ with $f(x_1) = f(x_2)$ define a
self-intersection, which is a condition on pairs of $0$-jets
(``multijets'').

The $1$-jet $j^1 f(x) = (x, f(x), Df(x))$ records both the value
and the derivative.  Rank drop ($\det Df = 0$) is a condition on the
$1$-jet.  The \emph{Thom--Boardman strata}
$\Sigma^i = \{j^1 f : \mathrm{rank}(Df) = \dim M - i\}$ organise
these singularities by their corank.

Higher jets ($r \geq 2$) capture curvature effects and higher-order
instabilities.

\paragraph{Multijets and self-intersections.}
To detect self-intersections, one needs to compare the map at two
different points simultaneously.  The \emph{multijet space}
$J^r_s(M,N)$ consists of $s$-tuples of $r$-jets at distinct points.
The \emph{multijet transversality theorem} asserts that for a generic
map, the multijet extension is transversal to any submanifold of the
multijet space---in particular, to the self-intersection diagonal.

\begin{theorem}[Multijet Transversality]
\label{thm:multijet}
For any submanifold $W \subset J^r(M,N)$, the set of smooth maps
$f$ with $j^r f \pitchfork W$ is residual.
\end{theorem}

\paragraph{Whitney stratifications.}
A \emph{Whitney stratification} of a set $X$ is a decomposition into
smooth manifolds (``strata'') $X = \bigsqcup X_i$, satisfying
regularity conditions (Whitney's conditions~(a) and~(b)) that
control how the strata fit together.  The degeneracy set of a feature
map admits such a stratification, and Thom's theorem extends to
stratified targets.

\paragraph{The Thom--Mather classification.}
The Thom--Mather theory classifies the generic singularities of smooth
maps $f : \R^m \to \R^n$ into a hierarchy---folds, cusps,
swallowtails, and so on---indexed by codimension.  For maps between
spaces of ``nice dimensions'' ($m$ and $n$ not too different), this
classification is exhaustive: every generic singularity is equivalent
to one of the normal forms.

\begin{keyidea}
Jets provide the language for classifying the types and severities of
singularity.  The jet transversality theorem says that generic maps
have the simplest possible singularities (those of lowest codimension),
and avoid the more severe ones.  In the statistical setting, this
translates to: for a generic kernel, the feature map has the mildest
possible degeneracies.
\end{keyidea}

\section{The Infinite-Dimensional Picture}
\label{sec:infinite_dim}

\textit{We discuss why the full distributional framework involves
  infinite-dimensional spaces, state the Sard--Smale and Abraham
  transversality theorems, and explain the role of the Fredholm
  condition as the mechanism that reduces infinite-dimensional
  transversality to finite-dimensional linear algebra.  This section
  expands~\cite{TransversalityPaper}, Section~6.}

\medskip

The finite-dimensional transversality theory of the preceding sections
applies when both the parameter space and the feature space are
finite-dimensional.  This covers the case of parametric statistical
models with finitely many weak moments.  But the full distributional
framework involves the kernel space $\cS(\R^d)$, which is
infinite-dimensional.

The classical transversality theory extends to the infinite-dimensional
setting through two fundamental results.

\begin{theorem}[Sard--Smale, 1965]
\label{thm:sard_smale}
Let $f : X \to Y$ be a $C^q$ Fredholm map between separable Banach
manifolds, with $q > \max(0, \mathrm{index}(f))$.  Then the set of
regular values of $f$ is residual in $Y$.
\end{theorem}

\begin{theorem}[Abraham, 1963]
\label{thm:abraham}
Let $F : X \times \Lambda \to Y$ be $C^q$ with $F \pitchfork S$.  If
each $F_\lambda$ is a Fredholm map, then $F_\lambda \pitchfork S$ for
a residual set of $\lambda \in \Lambda$.
\end{theorem}

The key concept is the \emph{Fredholm condition}.  A bounded linear
operator $L : X \to Y$ between Banach spaces is Fredholm if its
kernel $\ker L$ is finite-dimensional and its cokernel
$Y / \mathrm{Im}(L)$ is finite-dimensional.  The \emph{index} of $L$
is $\dim \ker L - \dim \mathrm{coker}(L)$.

\begin{keyidea}
The Fredholm condition is the mechanism by which infinite-dimensional
transversality reduces to finite-dimensional linear algebra.  A
Fredholm operator has finite-dimensional kernel and cokernel, so the
transversality question---whether the image of the derivative plus
the tangent space of $S$ spans the ambient space---becomes a
finite-rank calculation, even when the ambient spaces are
infinite-dimensional.
\end{keyidea}

In the parametric setting of these notes, the parameter space $\Theta$
is finite-dimensional ($\Theta \subset \R^p$), so each feature map
$\Phi_\varphi : \Theta \to \cF$ is automatically Fredholm.
Abraham's theorem then asserts: if the full map
$F(\theta, \varphi) = \Phi_\varphi(\theta)$ is transversal to a
degeneracy stratum $D$, then for a residual set of kernels,
$\Phi_\varphi \pitchfork D$.

Two important caveats must be noted.  First, the theorems of Smale and
Abraham are formulated for maps between separable \emph{Banach}
manifolds, whereas the Schwartz space $\cS(\R^d)$ carries a
Fr\'echet (not Banach) topology.  A full verification of the Fredholm
hypotheses in the distributional setting, or an adaptation to the
Fr\'echet framework, is a non-trivial analytic problem that is beyond
the scope of these notes.

Second, for semiparametric and nonparametric models---where the
parameter space itself is infinite-dimensional---the transversality
question becomes genuinely infinite-dimensional in both source and
target.  Quinn's extension of transversality to Banach manifolds
provides the appropriate framework, and the Fredholm condition
remains the key tool: it ensures that the obstruction to
transversality is always finite-dimensional, even when the ambient
spaces are not.  The development of these extensions is an active area
of research.

\begin{remark}[Statistical perspective]
\label{rem:infinite_stat}
The limitation to parametric models (finite-dimensional $\Theta$) is
the main statistical restriction of the present theory.  Extending to
semiparametric models---where $\theta = (\psi, \eta)$ with $\psi$
finite-dimensional and $\eta$ an infinite-dimensional nuisance
component---is the natural next step.  The Fredholm framework is
precisely designed for this: it reduces the infinite-dimensional
nuisance problem to a finite-rank calculation on the cokernel.
\end{remark}


\newpage
\part*{Part III: Statistical Theory Through the Transversality Lens}
\addcontentsline{toc}{part}{Part III: Statistical Theory Through the
  Transversality Lens}

\bigskip

\section{The Degeneracy Stratification}
\label{sec:degeneracy_strat}

\textit{We return to the statistical pathologies of Part~I and organise
  them into a geometric hierarchy: the degeneracy stratification of
  the feature map (cf.~\cite{TransversalityPaper}, Section~7).
  Each type of statistical pathology corresponds to a condition on
  the jets of the feature map, and the strata are ordered by
  codimension.}

\medskip

The pathologies discussed in Part~I---representation failure,
non-identifiability, singular information, moment indeterminacy, and
higher-order instabilities---can now be organised into a single
geometric framework: they are degeneracies of the feature map
$\Phi_\varphi$ and its jet extensions.

We identify five principal types of degeneracy.

\begin{center}
\begin{tabular}{llll}
\toprule
\textbf{Type} & \textbf{Degeneracy} & \textbf{Jet condition} &
  \textbf{Ref.} \\
\midrule
0 & Representation failure & Feature map undefined &
  \S\ref{sec:representation} \\
I & Non-identifiability & Self-intersection of $\Phi_\varphi$ &
  \S\ref{sec:identifiability} \\
II & Singular information & $\mathrm{rank}(D\Phi_\varphi) < p$ &
  \S\ref{sec:fisher_info} \\
III & Moment indeterminacy & Non-separation at distributional level &
  \S\ref{sec:identifiability} \\
IV & Higher-order instability & Degeneracies of $j^r\Phi_\varphi$,
  $r \geq 2$ & \S\ref{sec:jets} \\
\bottomrule
\end{tabular}
\end{center}

\medskip

\textbf{Type~0: Representation degeneracy.}  The feature map cannot
even be defined in the classical sense, because no density exists
(Section~\ref{sec:representation}).  In the distributional framework,
the kernel creates the feature map: the pairing
$\langle T_\theta, g_j \varphi \rangle$ is well-defined and smooth
even when no density is available.

\textbf{Type~I: Non-identifiability.}  The feature map has
self-intersections: $\Phi_\varphi(\theta_1) = \Phi_\varphi(\theta_2)$
with $\theta_1 \neq \theta_2$.  This is a condition on pairs of
$0$-jets in the multijet space
(Section~\ref{sec:jets}).

\textbf{Type~II: Singular information.}  The Jacobian
$D\Phi_\varphi(\theta)$ drops rank: the distributional metric tensor
$\wG(\theta)$ is singular.  This is a Thom--Boardman singularity
$\Sigma^1$ of the $1$-jet.

\textbf{Type~III: Moment indeterminacy.}  The feature map fails to
separate distributions (not just parameter values).  This corresponds
to M-indeterminacy in the classical moment problem: distinct
distributions are mapped to the same weak moment sequence.

\textbf{Type~IV: Higher-order instability.}  Conditions on higher
jets ($r \geq 2$): inflection points of weak moment functions,
vanishing curvature of the distributional metric, or instabilities
arising from complex dependency structures (such as non-chordal
graphical models).

These strata form a Whitney stratification of the degeneracy set, with
strata of increasing codimension:
$\mathrm{codim}(\text{Type~I}) < \mathrm{codim}(\text{Type~II})
< \mathrm{codim}(\text{Type~IV})$.  Transversality to a stratum of
codimension exceeding $\dim \Theta = p$ implies avoidance of that
stratum (Section~\ref{sec:transversality}).

\begin{keyidea}
The five types of statistical degeneracy form a hierarchy organised by
codimension.  Transversality theory predicts that for a generic kernel,
the feature map avoids all strata of sufficiently high codimension.
The ``sufficiently high'' threshold depends on the dimension of the
parameter space and the number of weak moments used.
\end{keyidea}

\section{The Main Theorem and What It Says}
\label{sec:main_theorem}

\textit{We state the finite-dimensional weak transversality theorem and
  explain its statistical content: for a generic kernel, the feature
  map avoids degeneracy strata of high codimension.  We spell out the
  codimension counting that determines how many weak moments are needed
  to guarantee identifiability and information regularity.}

\medskip

We can now state the main result, which instantiates the parametric
transversality theorem (Theorem~\ref{thm:parametric}) in the
distributional statistical setting
(cf.~\cite{TransversalityPaper}, Theorem~4.1 and Corollary~4.2).

\begin{theorem}[Finite-dimensional weak transversality]
\label{thm:main}
Let $\Theta \subset \R^p$ and $\Lambda \subset \R^q$ be open sets,
and let $\{T_\theta : \theta \in \Theta\} \subset \cS'(\R)$ be a
$C^r$ parametric distributional model.  Let
$\lambda \mapsto \varphi_\lambda \in \cS(\R)$ be a $C^r$
finite-dimensional family of positive Schwartz kernels.

Fix moment orders $0 \leq j_0 < \cdots < j_K$, and define the joint
weak moment feature map
\[
  F(\theta, \lambda) = \Phi_\lambda(\theta)
  = \bigl({}^{(\varphi_\lambda)}m_{j_k}(\theta)\bigr)_{k=0}^{K}.
\]
Assume that $F$ is $C^r$ and transversal to a smooth submanifold
$D \subset \R^{K+1}$.

Then for a residual subset $\Lambda_D \subset \Lambda$, the
restricted feature map $\Phi_\lambda$ is transversal to $D$ for
every $\lambda \in \Lambda_D$.  In particular, if
$\mathrm{codim}(D) > p$, then
$\Phi_\lambda(\Theta) \cap D = \varnothing$ for generic $\lambda$.
\end{theorem}

The assumption $F \pitchfork D$ can be verified in practice using
the differential criteria described in the next section
(corresponding to~\cite{TransversalityPaper}, Section~5), where it
reduces to rank conditions on the parameter and kernel derivatives
of the joint feature map.

\paragraph{What does the theorem say statistically?}

The theorem has two main consequences, obtained by applying the
codimension-counting principle of Section~\ref{sec:transversality}.

\textbf{Generic identifiability.}  The self-intersection diagonal in
the multijet space (Type~I degeneracy) has codimension $K+1$ in the
multijet space.  Thus if $K+1 > 2p$ (i.e.\ we use more than $2p$
weak moments), transversality implies that the feature map is
generically injective: distinct parameter values give rise to distinct
weak moment sequences.

\textbf{Generic information regularity.}  The rank-drop stratum
$\Sigma^1 = \{\mathrm{rank}(D\Phi) < p\}$ has codimension $K+1-p+1$
in the $1$-jet space.  Thus if $K+1 > 2p-1$, generically
$\det \wG(\theta) > 0$ for all $\theta$: the distributional
information matrix is non-singular everywhere.

\begin{keyidea}[Moment thresholds:]
For a $p$-parameter model, $K+1 \geq 2p$ weak moments generically
ensure identifiability, and $K+1 \geq 2p+1$ weak moments generically
ensure information regularity.  In practice, the weak characteristic
function (an infinite-dimensional feature) easily exceeds these
thresholds.
\end{keyidea}

\section{Verifiable Conditions}
\label{sec:verifiable}

\textit{We develop practical conditions---formulated as rank conditions
  on the Jacobian of the joint feature map---under which the
  transversality hypothesis of the main theorem can be verified.  We
  illustrate with four examples: the location family, the log-normal,
  Stein discrepancies, and graphical models.}

\medskip

The main theorem (Theorem~\ref{thm:main}) assumes that the joint map
$F(\theta, \lambda) = \Phi_\lambda(\theta)$ is transversal to the
degeneracy stratum $D$.  How can this assumption be checked in
practice?

The key observation is that the domain
$\Theta \times \Lambda$ is a product, so the derivative decomposes
into a \emph{model component} and a \emph{kernel component}:
\[
  DF(\theta, \lambda)
  = \bigl(D_\theta F, \;\; D_\lambda F\bigr).
\]
The model component $D_\theta F$ captures the intrinsic geometry of
the parametric family; the kernel component $D_\lambda F$ provides
supplementary directions from the kernel variation.

\begin{lemma}[Component-wise transversality criterion
  {\cite[Lemma~5.1]{TransversalityPaper}}]
\label{lem:component}
Let $\pi_N : \R^{K+1} \to N_y D$ denote the orthogonal projection
onto the normal space to $D$ at $y = F(\theta,\lambda) \in D$.  Then
$F \pitchfork D$ at $(\theta, \lambda)$ if and only if
\[
  \pi_N\bigl(\mathrm{Im}\, D_\theta F\bigr)
  + \pi_N\bigl(\mathrm{Im}\, D_\lambda F\bigr)
  = N_y D.
\]
\end{lemma}

The lemma says that transversality can be checked by examining the two
components separately: neither the model derivatives nor the kernel
derivatives need individually span the normal space, as long as their
normal projections together do.

The cleanest sufficient condition is \emph{submersivity}:

\begin{theorem}[Submersivity implies universal transversality
  {\cite[Theorem~5.3]{TransversalityPaper}}]
\label{thm:submersivity}
If the Jacobian $DF(\theta, \lambda)$ is surjective (rank $K+1$)
at every $(\theta, \lambda)$, then $F \pitchfork D$ for \emph{every}
smooth submanifold $D \subset \R^{K+1}$.
\end{theorem}

This is the strongest condition but also the easiest to check: one
only needs to verify that the $(K+1) \times (p+q)$ Jacobian matrix has
no rank deficiency, without knowing the specific stratum $D$.

When the model alone is degenerate, the kernel can compensate:

\begin{proposition}[Kernel-induced rank enrichment
  {\cite[Proposition~5.5]{TransversalityPaper}}]
\label{prop:kernel_rank}
If $D_\theta \Phi_\lambda(\theta)$ has rank $r < \min(p, K+1)$ and
$D_\lambda F$ contributes $\ell$ linearly independent directions
outside $\mathrm{Im}(D_\theta F)$, then
$\mathrm{rank}\, DF \geq r + \ell$.
\end{proposition}

\begin{keyidea}
The kernel acts as a source of supplementary directions that can lift
degeneracies the model alone cannot resolve.  The component-wise
criterion shows that the model and kernel contributions are
complementary: transversality holds whenever they jointly span the
normal space to the degeneracy stratum.
\end{keyidea}

\paragraph{Example 1: One-parameter location family.}
Consider $\{P_\mu : \mu \in \R\}$ with Gaussian kernel
$\varphi_s(x) = (2\pi s^2)^{-1/2} e^{-x^2/(2s^2)}$, $s > 0$.  The
joint map $F(\mu, s) = {}^{(\varphi_s)}w_0(\mu)$ has a $1 \times 2$
Jacobian $(\partial_\mu F, \partial_s F)$.  The $s$-derivative
$\partial_s F = \E_\mu[X^2 \varphi_s(X)/s^3]$ is strictly positive,
so $DF$ has rank $1$ everywhere.  Thus $F$ is a submersion and
Theorem~\ref{thm:submersivity} gives transversality to any degeneracy
stratum.

\paragraph{Example 2: The log-normal
  (cf.~\cite{TransversalityPaper}, Section~5.6.2 and
  Proposition~5.4).}
The two-parameter log-normal model with Gaussian kernel has a
$2 \times 3$ Jacobian (derivatives with respect to $\mu$, $\sigma$,
and $s$).  The $\mu$-derivative involves the centred log-score, which
is an odd function of $\ln X - \mu$; the $\sigma$-derivative involves
a quadratic, which is even.  The odd/even asymmetry ensures that the
$2 \times 2$ model Jacobian has rank $2$ (immersivity), and the
kernel derivative adds a third direction.  Thus the joint map is a
submersion and transversality holds universally.

\paragraph{Example 3: Stein discrepancies
  (cf.~\cite{TransversalityPaper}, Section~5.6.3 and
  Section~9).}
The weak Stein map requires two conditions for transversality: (a) the
test functions must be measure-determining (injectivity), and (b) the
joint Jacobian must be surjective.  Condition (a) alone guarantees
injectivity but not transversality; condition (b) is the additional
requirement.

\paragraph{Example 4: Gaussian graphical models
  (cf.~\cite{TransversalityPaper}, Section~5.6.4 and
  Example~8.4).}
For a Gaussian graphical model with precision matrix $\Omega$,
the model Jacobian with respect to the free entries of $\Omega$ has
full column rank for chordal graphs.  For non-chordal graphs (e.g.\ a
$4$-cycle), the model alone may be rank-deficient, but the kernel
derivatives restore full rank via
Proposition~\ref{prop:kernel_rank}.

\section{M-Indeterminacy Resolved}
\label{sec:M_indet_resolved}

\textit{We return to the log-normal example of
  Section~\ref{sec:identifiability} and show how the transversality
  framework explains the resolution of M-indeterminacy by the kernel.
  The Stieltjes perturbation is the tangent direction along which
  transversality fails; the Gaussian kernel destroys the symmetry
  responsible for the failure.  This section expands the treatment
  in~\cite{TransversalityPaper}, Section~8.1 and
  Proposition~8.1.}

\medskip

In Section~\ref{sec:identifiability} we saw that the log-normal
distribution is M-indeterminate: the classical moment map fails to
separate the Stieltjes class $\{\mu_a : a \in [-1,1]\}$ from the
log-normal $\mu$.  We can now understand this failure geometrically.

The classical moment map corresponds to the feature map with
$\varphi \equiv 1$ (no kernel).  The Stieltjes perturbation
$h(x) = \sin(2\pi \log x)$ is a tangent direction along which the
moment map is completely degenerate: all moments of $h$ with respect
to the log-normal vanish, so the perturbation is invisible to the
moment map.  In the language of Part~II, the moment map is
\emph{not transversal} to the self-intersection diagonal in the
multijet space.

The Gaussian kernel $\varphi(x) = (2\pi)^{-1/2}e^{-x^2/2}$ breaks
the symmetry responsible for this cancellation.  The integral
$\int x^n \sin(2\pi \log x)\, e^{-x^2/2}\, \mu(dx)$ no longer
vanishes, because the Gaussian weight $e^{-x^2/2}$ is not compatible
with the multiplicative periodicity of $\sin(2\pi \log x)$.  The
weak moment map separates the Stieltjes class from the log-normal, and
for generic kernels, the feature map is an immersion (the distributional
information matrix is non-singular).

\begin{keyidea}
The classical M-indeterminacy of the log-normal is a failure of
transversality.  The kernel restores transversality by breaking the
symmetry of the Stieltjes class.  The parametric transversality
theorem guarantees that this resolution is generic: it holds for
``almost every'' kernel, not just specific ones.
\end{keyidea}

\section{Information Regularity and Robustness}
\label{sec:info_robust}

\textit{We show that the distributional metric tensor is the first
  fundamental form of the feature-map immersion, connecting information
  regularity to immersivity and robustness to boundedness of the
  induced metric.}

\medskip

The distributional metric tensor
\[
  \wG_{ab}(\theta) = (D\Phi_\varphi)^\top (D\Phi_\varphi)
\]
is the first fundamental form of the immersion
$\Phi_\varphi : \Theta \to \cF$.  It measures infinitesimal distances
in the parameter space as induced by the feature map.

Transversality of the feature map to the rank-drop stratum $\Sigma^1$
ensures that $\wG(\theta)$ is non-singular for generic kernels
(Section~\ref{sec:main_theorem}).  This is the geometric content of
information regularity.

Robustness, in this framework, corresponds to \emph{boundedness} of
the metric tensor.  The kernel ensures that the weak moments
$\wm{j}(\theta) = \E_\theta[X^j \varphi(X)]$ involve the rapidly
decaying factor $\varphi(X)$, which bounds the contribution of extreme
observations.  As a consequence, the derivatives
$\partial \wm{j}/\partial \theta_a$ are bounded, and the metric tensor
$\wG_{ab}(\theta)$ is bounded.

\begin{keyidea}
Information regularity = the feature map is an immersion (the metric
tensor is non-degenerate).  Robustness = the metric tensor is bounded
(the manifold has finite geodesic lengths).  The kernel ensures both
properties simultaneously.
\end{keyidea}

The classical Fisher information geometry of Amari and
Barndorff-Nielsen is the special case $\varphi \equiv 1$; the
distributional information geometry is the general case
(see~\cite{TransversalityPaper}, Section~12.3).
In the classical limit, the metric may degenerate (singular
information) or diverge (unbounded influence); the kernel regularises
both pathologies.

\section{The Behrens--Fisher Problem as Nuisance Non-Transversality}
\label{sec:behrens_fisher_geom}

\textit{We return to the Behrens--Fisher problem of
  Section~\ref{sec:behrens_fisher_stat} and show that the impossibility
  of exact inference is a transversality failure: the null hypothesis
  is not in generic position relative to the nuisance-parameter
  fibration.  The kernel provides a family of deformations that
  restore transversality for generic kernel scale.  This section
  expands~\cite{TransversalityPaper}, Section~10.}

\medskip

The Behrens--Fisher problem has $\theta = (\mu_1, \mu_2, \sigma_1,
\sigma_2)$ and null hypothesis
$\Theta_0 = \{\mu_1 = \mu_2\} \cong \R \times (0,\infty)^2$.  The
nuisance parameters $(\sigma_1, \sigma_2)$ define a fibration
$\pi : \Theta \to \Psi = \{(\mu_1, \mu_2)\}$, and the difficulty is
that the sufficient-statistic projection is not transversal to
$\Theta_0$ relative to this fibration.

With a Gaussian kernel $\varphi_s$, the zeroth weak moments of the two
populations are
\[
  {}^{(\varphi_s)}w_0^{(k)}
  = \frac{1}{\sqrt{\sigma_k^2 + s^2}}\,
  \exp\!\left(-\frac{\mu_k^2}{2(\sigma_k^2 + s^2)}\right),
  \qquad k = 1, 2.
\]
The key feature is the coupling $\sigma_k^2 + s^2$: the location and
scale parameters are mixed through the kernel scale.  When
$s^2 \gg \max(\sigma_1^2, \sigma_2^2)$, the weak moments become
approximately
${}^{(\varphi_s)}w_0^{(k)} \approx s^{-1} e^{-\mu_k^2/(2s^2)}$, and
the nuisance parameters effectively disappear.

The kernel provides a one-parameter family of deformations of the
feature map.  By the parametric transversality theorem, for generic
$s$ the deformed null hypothesis is transversal to the nuisance
fibration.  In the classical limit $s \to \infty$, transversality is
lost: the ``paradox'' corresponds to a degenerate point in kernel space.

There is a trade-off: large $s$ gives nuisance insensitivity but
reduces statistical power (the feature map becomes coarse), analogous
to the efficiency--robustness trade-off.

\begin{keyidea}
The Behrens--Fisher problem is a non-transversality of the null
hypothesis relative to the nuisance-parameter fibration.  The kernel
resolves it by deforming the feature map into a generic position.
The classical framework corresponds to a degenerate point in the
space of representations.
\end{keyidea}

\section{Inferential Separation as Transversality}
\label{sec:inf_separation}

\textit{We show that the theory of inferential separation---sufficiency,
  ancillarity, and the nonformation principle---admits a natural
  transversality interpretation.  The Bhapkar--Godambe projection
  enforces transversality to the nuisance tangent space, and
  sinusoidal inference functions achieve transversality automatically.
  This section expands~\cite{TransversalityPaper}, Section~12.7.}

\medskip

Consider a parametric model with $\theta = (\psi, \phi)$, where
$\psi$ is the interest parameter and $\phi$ is a nuisance parameter.
The nuisance parameter defines a fibration
$\pi : \Theta \to \Psi$, $\pi(\psi, \phi) = \psi$, whose fibres
$\{\psi\} \times \Xi$ are the nuisance orbits.

Inferential separation requires that inference about $\psi$ be
insensitive to $\phi$.  In the feature-map language, this means that
the restriction of $\Phi_\varphi$ to the interest direction is
transversal to the nuisance fibres: the image of the interest subspace
intersects the nuisance tangent space only trivially.

In the inference-function formulation, an inference function for
$\psi$ lives in the orthogonal complement $\mathcal{T}_N^\perp$ of
the nuisance tangent space.  The \emph{Bhapkar--Godambe projection}
takes an arbitrary quasi-inference function and projects it onto
$\mathcal{T}_N^\perp$---this is precisely the operation of deforming
the inference function into a transversal position.

A remarkable fact in the distributional framework is that
\emph{sinusoidal inference functions}
$\psi_c(x, \mu) = \sin(c(x - \mu))$ for symmetric location-scale
models satisfy the orthogonality condition \emph{automatically},
without requiring explicit projection.  This automatic transversality
arises from the symmetry of the characteristic function.

The hierarchy of nonformation concepts corresponds to progressively
weaker transversality conditions:

\textbf{S-nonformation:} the feature map factors completely through
the interest projection (full transversality to nuisance fibres).

\textbf{I-nonformation:} the conditional feature map is saturated
(local transversality at the maximum).

\textbf{L-nonformation:} the profile depends on the data only through
a reduction (transversality at the level of the normed profile).

\begin{keyidea}
The various notions of inferential separation are manifestations of a
single geometric principle: the feature map being in generic position
relative to the nuisance structure.  The kernel provides a mechanism
for achieving transversality even when classical likelihood-based
separation fails.
\end{keyidea}

\section{The Singular Limit and the Classical Framework}
\label{sec:singular_limit}

\textit{We interpret the classical statistical framework as a singular
  (degenerate) limit of the distributional framework: the limit
  $\varphi_s \to 1$ is a path in kernel space that leaves the residual
  set of transversal kernels and enters a degeneracy stratum.  The
  classical pathologies are not fundamental but arise from working at
  a degenerate point
  (cf.~\cite{TransversalityPaper}, Section~7.4).}

\medskip

Throughout these notes, we have seen that the classical framework
(no kernel, or equivalently $\varphi \equiv 1$) is the setting in
which pathologies arise: M-indeterminacy, singular information,
non-robustness, the Behrens--Fisher paradox.  The distributional
framework, with a non-trivial kernel, resolves each of these.

The transversality perspective provides a precise explanation.  The
kernel space $\cS(\R^d)$ contains a residual set of ``good'' kernels
for which the feature map is transversal to all degeneracy strata.
The constant function $\varphi \equiv 1$ does not belong to
$\cS(\R^d)$ (it is not rapidly decaying), but it can be approximated
by a family $\varphi_s$ with $\varphi_s \to 1$ as $s \to \infty$.

This limit is a path in (an extension of) the kernel space that
\emph{exits} the residual set of transversal kernels and
\emph{enters} a degeneracy stratum.  The classical framework is not
the ``natural'' setting; it is a \emph{degenerate point} in the space
of representations.

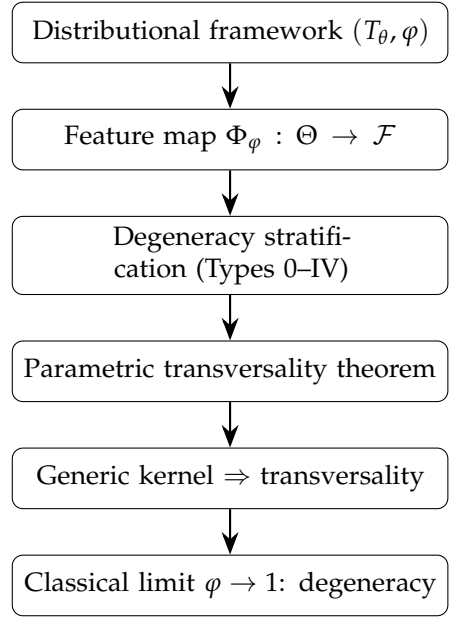
\begin{keyidea}[Conceptual chain:]
\begin{center}
\begin{tikzpicture}[
  node distance=0.6cm,
  box/.style={draw, rounded corners, minimum height=0.8cm,
    text width=5.5cm, align=center, font=\small},
  arr/.style={-{Stealth[length=3mm]}, thick}
]
\node[box] (A) {Distributional framework $(T_\theta, \varphi)$};
\node[box, below=of A] (B) {Feature map $\Phi_\varphi : \Theta \to \cF$};
\node[box, below=of B] (C) {Degeneracy stratification (Types 0--IV)};
\node[box, below=of C] (D) {Parametric transversality theorem};
\node[box, below=of D] (E) {Generic kernel $\Rightarrow$ transversality};
\node[box, below=of E] (F) {Classical limit $\varphi \to 1$: degeneracy};
\draw[arr] (A) -- (B);
\draw[arr] (B) -- (C);
\draw[arr] (C) -- (D);
\draw[arr] (D) -- (E);
\draw[arr] (E) -- (F);
\end{tikzpicture}
\end{center}
\end{keyidea}

This picture is analogous to a resolution of singularities in
algebraic geometry: the distributional framework provides a
one-parameter deformation (indexed by the kernel) that resolves the
degeneracies of the classical framework.  The classical results are
recovered as limiting cases, and the pathologies are understood as
artefacts of working at the degenerate point $\varphi \equiv 1$.

\begin{remark}[Statistical perspective]
\label{rem:singular_stat}
The singular-limit perspective suggests a practical principle: when a
classical statistical method encounters difficulties (non-robustness,
moment problems, nuisance effects), one should consider whether these
difficulties arise from the degenerate nature of the classical
representation, and whether a kernel regularisation can resolve them.
\end{remark}

\begin{remark}[Geometric perspective]
\label{rem:singular_geom}
The degenerate limit $\varphi \to 1$ is not a pathology of the
distributional framework; it is a \emph{confirmation} of the
transversality picture.  The theory predicts that the set of
non-transversal kernels is meagre (a countable union of nowhere dense
sets), and the classical ``kernel'' $\varphi \equiv 1$ is a specific
non-transversal point.  The fact that classical statistics encounters
difficulties at this point is exactly what the theory predicts.
\end{remark}

\section{Conceptual Summary}
\label{sec:conceptual_summary}

\textit{We summarise the main ideas of the notes, emphasising the
  conceptual chain linking statistical pathologies, the feature map,
  transversality, and the role of the kernel.}

\medskip

The starting point of these notes was a collection of classical
statistical pathologies: models without densities, moment
indeterminacy, singular information, non-robustness, and the
Behrens--Fisher problem.  These are usually treated as separate
difficulties, each requiring its own ad hoc solution.

The distributional framework, based on distribution--kernel pairs
$(T_\theta, \varphi)$, provides a unified setting in which all these
pathologies are resolved.  The kernel induces a feature map
$\Phi_\varphi : \Theta \to \cF$, and the statistical properties of
the model are encoded in the geometry of this map.

The transversality perspective explains \emph{why} the kernel works:
it acts as a generic perturbation that places the feature map in a
transversal position relative to the degeneracy strata.  The
parametric transversality theorem guarantees that this resolution is
generic: for ``almost every'' kernel, the feature map avoids the
degeneracies.

The classical framework corresponds to the degenerate limit
$\varphi \to 1$, where transversality is lost.  The statistical
pathologies are not fundamental; they arise from working at a
degenerate point in the space of representations.

\begin{keyidea}[Summary of the conceptual chain:]
\begin{center}
\begin{tabular}{rcl}
Statistical pathology & $\longleftrightarrow$ &
  geometric degeneracy of $\Phi_\varphi$ \\[4pt]
Kernel regularisation & $\longleftrightarrow$ &
  generic perturbation \\[4pt]
Transversality theorem & $\longleftrightarrow$ &
  genericity of resolution \\[4pt]
Classical framework & $\longleftrightarrow$ &
  degenerate limit
\end{tabular}
\end{center}
\end{keyidea}

\newpage
\section{Exercises}
\label{sec:exercises}

\begin{exercise}[Sard's theorem and the sphere]
Let $f : \R^3 \to \R$ be defined by $f(x,y,z) = x^2 + y^2 + z^2$.
Compute the set of critical values of $f$ and verify that it has
measure zero.  For each regular value $r > 0$, describe $f^{-1}(r)$
as a smooth manifold.
\end{exercise}

\begin{exercise}[Transversality of a line and a surface]
Let $\gamma : \R \to \R^3$ be the line $\gamma(t) = (t, 0, 1)$ and
let $S = \{(x,y,z) \in \R^3 : z = 0\}$ be the $xy$-plane.  Is
$\gamma \pitchfork S$?  What if
$\gamma(t) = (t, 0, t^2)$?  Discuss.
\end{exercise}

\begin{exercise}[Codimension counting]
A parametric model has $p = 3$ parameters and uses $K+1 = 8$ weak
moments.  What is the expected codimension of the non-identifiability
stratum in the multijet space?  Is $K+1$ large enough to guarantee
generic identifiability?  What about generic information regularity?
\end{exercise}

\begin{exercise}[M-indeterminacy and transversality]
Verify that $\int_0^\infty x^n \sin(2\pi \log x)\, \mu(dx) = 0$ for
every $n \geq 0$, where $\mu$ is the standard log-normal distribution,
by writing $x^n = \exp(n \log x)$ and computing the resulting Gaussian
integral.  Explain why this cancellation fails when the integrand
includes a Gaussian kernel $\varphi(x) = e^{-x^2/2}$.
\end{exercise}

\begin{exercise}[One-parameter location family]
Consider the Cauchy location family
$f(x;\mu) = \frac{1}{\pi}\frac{1}{1+(x-\mu)^2}$ with Gaussian
kernel $\varphi_s(x) = (2\pi s^2)^{-1/2} e^{-x^2/(2s^2)}$.  Compute
the zeroth weak moment ${}^{(\varphi_s)}w_0(\mu)$ and verify that the
joint Jacobian $DF(\mu, s)$ has rank $1$ at every point.
\end{exercise}

\begin{exercise}[Behrens--Fisher regularisation]
For the Behrens--Fisher problem with Gaussian kernel, verify the
formula
\[
  {}^{(\varphi_s)}w_0^{(k)}
  = \frac{1}{\sqrt{\sigma_k^2 + s^2}}\,
  \exp\!\left(-\frac{\mu_k^2}{2(\sigma_k^2 + s^2)}\right).
\]
Show that for $s^2 \gg \sigma_k^2$, the dependence on $\sigma_k^2$
becomes negligible.  Discuss the trade-off between nuisance
insensitivity and statistical power.
\end{exercise}

\begin{exercise}[Submersivity]
Let $F : \R^2 \times \R \to \R^2$ be a smooth map with
$DF(\theta, \lambda)$ a $2 \times 3$ matrix.  Give a necessary and
sufficient condition on $DF$ for $F$ to be a submersion.  Show that
if $F$ is a submersion, then $F \pitchfork D$ for every smooth
submanifold $D \subset \R^2$.
\end{exercise}

\begin{exercise}[Fredholm operators]
Let $L : \ell^2 \to \ell^2$ be the right shift operator
$L(x_1, x_2, \ldots) = (0, x_1, x_2, \ldots)$.  Compute $\ker L$
and $\mathrm{coker}(L)$.  Is $L$ Fredholm?  What is its index?
Compare with the left shift $R(x_1, x_2, \ldots) = (x_2, x_3, \ldots)$.
\end{exercise}


\newpage
\appendix

\section{Review of Geometric Prerequisites}
\label{app:prerequisites}

\textit{This appendix collects the basic notions from differential
  geometry and topology that are used throughout these notes.  It is
  intended for readers with a background in analysis and linear
  algebra who have not had a course in differential topology.  For
  comprehensive treatments, see
  Guillemin--Pollack~\cite{GuilleminPollack1974} or
  Hirsch~\cite{Hirsch1976}.}

\subsection{Smooth manifolds and submanifolds}
\label{app:manifolds}

A \emph{smooth manifold} of dimension~$m$ is, informally, a space
that locally looks like $\R^m$.  More precisely, it is a topological
space $M$ equipped with a collection of homeomorphisms
$\phi_\alpha : U_\alpha \to V_\alpha$ (called \emph{charts}), where
$U_\alpha \subset M$ are open sets covering $M$ and
$V_\alpha \subset \R^m$ are open, such that the \emph{transition
  maps} $\phi_\beta \circ \phi_\alpha^{-1}$ are smooth ($C^\infty$)
wherever defined.

\paragraph{Examples.}
(i)~Every open subset $U \subset \R^m$ is a smooth manifold of
dimension~$m$ (with the identity as its single chart).  In particular,
the parameter space $\Theta \subset \R^p$ of a statistical model is a
smooth manifold.
(ii)~The $n$-sphere
$S^n = \{x \in \R^{n+1} : \|x\| = 1\}$ is a smooth manifold of
dimension~$n$, covered by two charts (stereographic projections from
the north and south poles).
(iii)~More generally, if $f : \R^m \to \R^k$ is smooth and
$y$ is a \emph{regular value} (see Section~\ref{sec:smooth_maps}),
then $f^{-1}(y)$ is a smooth manifold of dimension~$m - k$.  This is
how many manifolds arise in practice.

A \emph{smooth submanifold} of $M$ of dimension $d$ is a subset
$S \subset M$ that is itself a smooth manifold of dimension~$d$, and
whose inclusion $S \hookrightarrow M$ is a smooth map with injective
derivative at every point.  The \emph{codimension} of $S$ in $M$ is
$\mathrm{codim}(S) = \dim M - \dim S$.

\paragraph{Why codimension matters.}
Two submanifolds of complementary dimension ``generically'' intersect
in isolated points.  If the codimension of $S$ exceeds $\dim M$, then
$S$ is empty (since $\dim S < 0$ is impossible).  This simple
dimension-counting principle is the basis of the transversality
arguments in these notes: if a degeneracy locus has high enough
codimension, a generic map avoids it entirely.

\subsection{Tangent spaces and the derivative of a smooth map}
\label{app:tangent}

The \emph{tangent space} $T_x M$ to a smooth manifold $M$ at a point
$x \in M$ is the vector space of ``velocity vectors'' of smooth
curves passing through~$x$.  Formally, two smooth curves
$\gamma_1, \gamma_2 : (-\varepsilon, \varepsilon) \to M$ with
$\gamma_1(0) = \gamma_2(0) = x$ are \emph{equivalent} if
$(\phi \circ \gamma_1)'(0) = (\phi \circ \gamma_2)'(0)$ in some (and
hence every) chart $\phi$ around $x$.  The tangent space $T_x M$ is
the set of equivalence classes.

If $M$ is an open subset of $\R^m$, then $T_x M \cong \R^m$
canonically: the tangent vector of $\gamma$ at $x$ is simply
$\gamma'(0) \in \R^m$.  This is the case most relevant to these
notes, since parameter spaces $\Theta \subset \R^p$ are open subsets
of Euclidean space.

If $M$ is a submanifold of $\R^n$ (for instance, a level set
$f^{-1}(y)$), then $T_x M$ is a \emph{linear subspace} of $\R^n$:
the set of all velocity vectors of curves in $M$ through $x$.  For
example, the tangent space to the sphere $S^2$ at a point $x$ is the
plane perpendicular to $x$ in $\R^3$.

\paragraph{The derivative of a smooth map.}
Given a smooth map $f : M \to N$ and a point $x \in M$, the
\emph{derivative} (or \emph{differential}) of $f$ at $x$ is the
linear map
\[
  Df_x : T_x M \to T_{f(x)} N
\]
defined by $Df_x([\gamma]) = [f \circ \gamma]$: it sends the velocity
vector of a curve $\gamma$ through $x$ to the velocity vector of the
image curve $f \circ \gamma$ through $f(x)$.

When $M \subset \R^m$ and $N \subset \R^n$ are open sets, $Df_x$ is
simply the \emph{Jacobian matrix}: the $n \times m$ matrix of partial
derivatives $(\partial f_i / \partial x_j)_{i,j}$.  This is the
concrete object that appears throughout these notes.

\paragraph{The transversality condition, revisited.}
With this language, the transversality condition
(Definition~\ref{def:transversality}) becomes concrete: $f \pitchfork
S$ at $x \in f^{-1}(S)$ means that the image of the Jacobian matrix
$Df_x$, together with the tangent space $T_{f(x)} S$ (a linear
subspace of $\R^n$), spans the full space $\R^n$.  This is a
\emph{rank condition} on a specific matrix, which can be checked by
linear algebra.

\subsection{Rank, immersions, and submersions}
\label{app:rank}

Let $L : \R^m \to \R^n$ be a linear map (i.e.\ a matrix).  The
\emph{rank} of $L$ is the dimension of its image:
$\mathrm{rank}(L) = \dim \mathrm{Im}(L)$.  By elementary linear
algebra, $\mathrm{rank}(L) \leq \min(m, n)$.

\begin{definition}
A linear map $L : \R^m \to \R^n$ is:
\begin{itemize}
\item \emph{injective} (or \emph{one-to-one}) if $\ker L = \{0\}$,
  equivalently $\mathrm{rank}(L) = m$;
\item \emph{surjective} (or \emph{onto}) if $\mathrm{Im}(L) = \R^n$,
  equivalently $\mathrm{rank}(L) = n$.
\end{itemize}
\end{definition}

These notions extend to smooth maps via the derivative:

\begin{definition}
A smooth map $f : M \to N$ is:
\begin{itemize}
\item an \emph{immersion} at $x$ if $Df_x$ is injective
  ($\mathrm{rank}(Df_x) = \dim M$);
\item a \emph{submersion} at $x$ if $Df_x$ is surjective
  ($\mathrm{rank}(Df_x) = \dim N$).
\end{itemize}
If the condition holds at every point, $f$ is called an immersion
(resp.\ submersion) globally.
\end{definition}

\paragraph{Statistical interpretation.}
In these notes, the feature map
$\Phi_\varphi : \Theta \to \R^{K+1}$ maps the $p$-dimensional
parameter space into the $(K+1)$-dimensional feature space.

\emph{Immersivity} ($\mathrm{rank}(D\Phi_\varphi) = p$ everywhere)
means that infinitesimal changes in the parameter always produce
detectable changes in the weak moments.  This is \emph{information
  regularity}: the distributional metric tensor
$\wG(\theta) = (D\Phi_\varphi)^\top (D\Phi_\varphi)$ is
non-singular.

\emph{Submersivity} of the joint map
$F(\theta, \lambda) = \Phi_\lambda(\theta)$ means that the combined
model-and-kernel Jacobian has rank $K+1$.  By
Theorem~\ref{thm:submersivity}, this implies transversality to
\emph{every} degeneracy stratum simultaneously.

\subsection{Residual sets and the Baire category theorem}
\label{app:residual}

The transversality theorems assert that ``generic'' maps are
transversal.  The precise meaning of ``generic'' is captured by the
topological notion of a \emph{residual set}.

\begin{definition}
Let $X$ be a topological space.  A subset $A \subset X$ is:
\begin{itemize}
\item \emph{nowhere dense} if its closure has empty interior:
  $\mathrm{int}(\overline{A}) = \emptyset$;
\item \emph{meagre} (or \emph{of first category}) if it is a
  countable union of nowhere dense sets;
\item \emph{residual} (or \emph{comeagre}) if its complement is
  meagre, equivalently, if it is a countable intersection of open
  dense sets.
\end{itemize}
\end{definition}

\paragraph{Intuition.}
Meagre sets are the topological analogue of measure-zero sets:
they are ``negligible'' in a topological sense.  Residual sets are
``large'': they contain ``almost all'' points.  However, unlike
measure-zero sets, meagre sets can have full measure, and residual
sets can have measure zero (for example, the irrational numbers form a
residual set in $\R$).

\begin{theorem}[Baire category theorem]
In a complete metric space (or more generally, a locally compact
Hausdorff space), every residual set is dense.  In particular, the
intersection of countably many open dense sets is dense.
\end{theorem}

\paragraph{Why this matters.}
The space $C^\infty(M, N)$ with the Whitney topology is a Baire space
(countable intersections of open dense sets are dense).  Thom's
transversality theorem asserts that the set of maps transversal to a
given submanifold is residual---hence dense---in this space.
Similarly, in the parametric setting, the set of kernels $\lambda$ for
which $\Phi_\lambda \pitchfork D$ is residual in the parameter space
$\Lambda$.  When $\Lambda$ is an open subset of $\R^q$ (which is
complete), the Baire category theorem guarantees that this residual
set is dense: transversal kernels are everywhere dense among all
kernels.

\subsection{Banach and Fr\'echet spaces}
\label{app:banach_frechet}

The infinite-dimensional extensions of transversality theory
(Section~\ref{sec:infinite_dim}) involve function spaces that are not
finite-dimensional.  Two classes of such spaces arise naturally.

A \emph{Banach space} is a complete normed vector space: a vector
space $X$ equipped with a norm $\|\cdot\|$ such that every Cauchy
sequence converges.  Examples include $\R^n$ (with any norm),
$L^p$ spaces, and $\ell^p$ spaces.

A \emph{Fr\'echet space} is a complete topological vector space whose
topology is defined by a \emph{countable} family of seminorms (rather
than a single norm).  The Schwartz space $\cS(\R^d)$ of rapidly
decaying smooth functions is a Fr\'echet space: its topology is
defined by the seminorms
$\|f\|_{\alpha,\beta} = \sup_{x} |x^\alpha D^\beta f(x)|$, indexed
by multi-indices $\alpha, \beta$.  It is \emph{not} a Banach space,
because no single norm captures the topology.

\paragraph{Why the distinction matters.}
The Sard--Smale theorem and Abraham's transversality theorem
(Theorems~\ref{thm:sard_smale} and~\ref{thm:abraham}) are formulated
for maps between separable \emph{Banach} manifolds.  The kernel space
$\cS(\R^d)$ is Fr\'echet, not Banach.  In the parametric setting of
these notes, this subtlety does not arise: the parameter space
$\Theta \subset \R^p$ is finite-dimensional, and the feature map
$\Phi_\varphi$ maps $\Theta$ into $\R^{K+1}$.  The
infinite-dimensional nature of $\cS(\R^d)$ enters only when we
consider the full kernel space as a parameter space in the parametric
transversality theorem, and in this context one works with
finite-dimensional subfamilies
$\Lambda \subset \R^q \hookrightarrow \cS(\R^d)$, for which the
Banach-space machinery is not needed.

\paragraph{Fredholm operators.}
A bounded linear operator $L : X \to Y$ between Banach spaces is
\emph{Fredholm} if $\ker L$ is finite-dimensional and
$Y / \mathrm{Im}(L)$ (the cokernel) is finite-dimensional.  The
\emph{Fredholm index} is
$\mathrm{index}(L) = \dim \ker L - \dim \mathrm{coker}(L)$.

The importance of the Fredholm condition is that it reduces
infinite-dimensional problems to finite-dimensional linear algebra:
the transversality question (``does the image of $Df$ plus the
tangent space of $S$ span the ambient space?'') involves only the
finite-dimensional kernel and cokernel.  This is why
Abraham's theorem requires the maps $F_\lambda$ to be Fredholm:
the condition ensures that the obstruction to transversality is
finite-dimensional, even when the ambient spaces are not.
See Exercise~20.8 and its solution (Appendix~\ref{app:solutions})
for concrete examples.

\newpage
\section{Solutions to Exercises}
\label{app:solutions}

\paragraph{Solution to Exercise~20.1 (Sard's theorem and the sphere).}
The Jacobian of $f(x,y,z) = x^2 + y^2 + z^2$ is
$Df(x,y,z) = (2x,\; 2y,\; 2z)$, which is a $1 \times 3$ matrix.
This is surjective (has rank~$1$) whenever $(x,y,z) \neq (0,0,0)$.
Hence the only critical point is the origin, and the only critical
value is $f(0,0,0) = 0$.  The set of critical values is $\{0\}$,
which has Lebesgue measure zero in $\R$, confirming Sard's theorem.

For each $r > 0$, the value $r$ is regular, and by the preimage
theorem (Theorem~\ref{thm:preimage}), $f^{-1}(r)$ is a smooth
submanifold of $\R^3$ of dimension $3 - 1 = 2$.  Concretely,
$f^{-1}(r) = \{(x,y,z) : x^2 + y^2 + z^2 = r\}$ is the $2$-sphere
of radius~$\sqrt{r}$.

\paragraph{Solution to Exercise~20.2 (Transversality of a line and
  a surface).}
The surface is $S = \{z = 0\}$, with tangent space
$T_p S = \mathrm{span}\{e_1, e_2\}$ at every point.  The ambient
space is $\R^3$.

\emph{Case 1:} $\gamma(t) = (t, 0, 1)$.  Then
$\gamma(\R) \cap S = \emptyset$ (the line lies entirely in the plane
$z = 1$).  Since $\gamma^{-1}(S)$ is empty, $\gamma \pitchfork S$
holds vacuously.

\emph{Case 2:} $\gamma(t) = (t, 0, t^2)$.  The intersection occurs
at $t = 0$, where $\gamma(0) = (0,0,0) \in S$.  The velocity vector
is $\gamma'(t) = (1, 0, 2t)$, so $\gamma'(0) = (1, 0, 0)$.  The
transversality condition requires
$\mathrm{span}\{\gamma'(0)\} + T_{\gamma(0)} S = \R^3$, i.e.\
$\mathrm{span}\{(1,0,0)\} + \mathrm{span}\{e_1, e_2\} = \R^3$.  But
$\mathrm{span}\{(1,0,0), e_1, e_2\} = \mathrm{span}\{e_1, e_2\}$,
which does not contain $e_3$.  Hence $\gamma \not\pitchfork S$: the
curve is tangent to $S$ at the intersection (it touches the plane
quadratically rather than crossing it transversally).

\paragraph{Solution to Exercise~20.3 (Codimension counting).}
With $p = 3$ and $K + 1 = 8$:

\emph{Identifiability:} The self-intersection diagonal in the
multijet space $J^0_2(\Theta, \R^{K+1})$ has codimension $K + 1 = 8$.
Since $8 > 2p = 6$, transversality implies that the feature map is
generically injective.  So $K + 1 = 8$ is sufficient for generic
identifiability.

\emph{Information regularity:} The Thom--Boardman stratum $\Sigma^1$
(rank drop by $1$) has codimension $K + 1 - p + 1 = 8 - 3 + 1 = 6$.
Since $6 > p = 3$, generically $\det \wG(\theta) > 0$ everywhere.
Alternatively, the condition $K + 1 > 2p - 1 = 5$ is satisfied
($8 > 5$), confirming generic information regularity.

\paragraph{Solution to Exercise~20.4 (M-indeterminacy and
  transversality).}
Let $\mu$ be the standard log-normal distribution, so
$\mu(dx) = (2\pi x^2)^{-1/2}\exp(-(\log x)^2/2)\,dx$ for $x > 0$.
We need to show that
$I_n = \int_0^\infty x^n \sin(2\pi \log x)\, \mu(dx) = 0$.

Substituting $y = \log x$ (so $x = e^y$, $dx = e^y \, dy$),
\begin{align*}
  I_n
  &= \int_{-\infty}^{\infty} e^{ny} \sin(2\pi y)\,
     \frac{1}{\sqrt{2\pi}} e^{-y^2/2}\, dy \\
  &= \frac{1}{\sqrt{2\pi}} \int_{-\infty}^{\infty}
     \sin(2\pi y)\, e^{ny - y^2/2}\, dy \\
  &= \frac{1}{\sqrt{2\pi}} \,\mathrm{Im}\!\int_{-\infty}^{\infty}
     e^{2\pi i y}\, e^{ny - y^2/2}\, dy.
\end{align*}
Completing the square: $ny - y^2/2 = -(y - n)^2/2 + n^2/2$, so
\[
  I_n = \frac{e^{n^2/2}}{\sqrt{2\pi}} \,\mathrm{Im}\!
  \int_{-\infty}^{\infty}
  e^{2\pi i y}\, e^{-(y-n)^2/2}\, dy
  = e^{n^2/2}\, \mathrm{Im}\bigl(e^{2\pi i n - 2\pi^2}\bigr)
  = e^{n^2/2 - 2\pi^2}\, \sin(2\pi n) = 0,
\]
since $\sin(2\pi n) = 0$ for every integer $n$.

With the Gaussian kernel $\varphi(x) = e^{-x^2/2}$, the analogous
integral becomes
$J_n = \int_0^\infty x^n \sin(2\pi \log x)\, e^{-x^2/2}\, \mu(dx)$.
After substituting $y = \log x$, the exponent in the integrand
contains $-e^{2y}/2$ (from $\varphi$), which is not a quadratic in
$y$.  The completing-the-square argument that produced
$\sin(2\pi n) = 0$ no longer applies: the Fourier analysis that
exploited the multiplicative periodicity of $\sin(2\pi \log x)$
relative to the log-normal is destroyed by the additive Gaussian
factor $e^{-x^2/2}$.  Hence $J_n \neq 0$ in general, and the kernel
breaks the M-indeterminacy.

\paragraph{Solution to Exercise~20.5 (One-parameter location family).}
The Cauchy location family has density
$f(x;\mu) = \frac{1}{\pi}\frac{1}{1+(x-\mu)^2}$.  The zeroth weak
moment with Gaussian kernel $\varphi_s$ is
\begin{align*}
  {}^{(\varphi_s)}w_0(\mu)
  &= \int_{-\infty}^{\infty}
     \frac{1}{\pi}\frac{1}{1+(x-\mu)^2}\cdot
     \frac{1}{\sqrt{2\pi s^2}}\, e^{-x^2/(2s^2)}\, dx.
\end{align*}
This is the convolution of a Cauchy density (centred at $\mu$) with a
Gaussian, evaluated at $0$.  By the known convolution formula,
$\text{Cauchy}(\mu, 1) * \text{Gaussian}(0, s^2)$ has the Voigt
profile, and the result can be expressed as the real part of a scaled
complementary error function.  What matters here is that the result is
a smooth function of $(\mu, s)$.

For the rank check: $\partial_\mu F = \partial_\mu\, {}^{(\varphi_s)}w_0(\mu)$
involves the odd part of the integrand (differentiating the Cauchy
density shifts mass), and $\partial_s F$ involves $\E_\mu[X^2
\varphi_s(X)/s^3] > 0$ (a moment of a positive function).  Since
$\partial_s F > 0$ everywhere, the $1 \times 2$ Jacobian
$DF(\mu, s) = (\partial_\mu F, \partial_s F)$ has rank~$1$ at every
point.

\paragraph{Solution to Exercise~20.6 (Behrens--Fisher regularisation).}
For population $k$ with $X \sim N(\mu_k, \sigma_k^2)$ and kernel
$\varphi_s(x) = (2\pi s^2)^{-1/2} e^{-x^2/(2s^2)}$,
\begin{align*}
  {}^{(\varphi_s)}w_0^{(k)}
  &= \int_{-\infty}^{\infty}
     \frac{1}{\sqrt{2\pi\sigma_k^2}}\, e^{-(x-\mu_k)^2/(2\sigma_k^2)}
     \cdot
     \frac{1}{\sqrt{2\pi s^2}}\, e^{-x^2/(2s^2)}\, dx.
\end{align*}
The integrand is the product of two Gaussians.  Combining the
exponents:
\[
  -\frac{(x-\mu_k)^2}{2\sigma_k^2} - \frac{x^2}{2s^2}
  = -\frac{x^2(\sigma_k^2 + s^2) - 2x\mu_k s^2 + \mu_k^2 s^2}
         {2\sigma_k^2 s^2}.
\]
Completing the square in $x$ and integrating the resulting Gaussian,
we obtain
\[
  {}^{(\varphi_s)}w_0^{(k)}
  = \frac{1}{\sqrt{2\pi(\sigma_k^2 + s^2)}}\,
  \exp\!\left(-\frac{\mu_k^2}{2(\sigma_k^2 + s^2)}\right).
\]
(The prefactor $1/\sqrt{2\pi}$ cancels with the normalisation of the
integrated Gaussian.)

When $s^2 \gg \sigma_k^2$, we have
$\sigma_k^2 + s^2 \approx s^2$, so
${}^{(\varphi_s)}w_0^{(k)} \approx (2\pi s^2)^{-1/2}\,
e^{-\mu_k^2/(2s^2)}$, which depends only on $\mu_k$ and $s$, not on
$\sigma_k$.  The nuisance parameter $\sigma_k^2$ disappears from the
leading-order expression.

\emph{Trade-off:} As $s$ increases, the feature map becomes
insensitive to $\sigma_k^2$ (good for nuisance elimination) but also
less sensitive to differences in $\mu_k$ (the exponent
$\mu_k^2/(2s^2)$ shrinks).  For large $s$,
${}^{(\varphi_s)}w_0^{(k)} \approx (2\pi s^2)^{-1/2}$ for all
parameter values: the feature map becomes nearly constant, and
statistical power vanishes.  This is the efficiency--robustness
trade-off discussed in~\cite{B}.

\paragraph{Solution to Exercise~20.7 (Submersivity).}
A smooth map $F : \R^2 \times \R \to \R^2$ has a $2 \times 3$
Jacobian $DF(\theta, \lambda)$ at each point.  The map is a
submersion if and only if $\mathrm{rank}\, DF(\theta, \lambda) = 2$ at
every $(\theta, \lambda)$, i.e.\ the $2 \times 3$ matrix has full
row rank.  Equivalently, at least one $2 \times 2$ minor of $DF$ is
non-zero at every point.

If $F$ is a submersion, then $\mathrm{Im}\, DF(\theta, \lambda)
= \R^2$ at every point.  For any smooth submanifold
$D \subset \R^2$, the transversality condition
$\mathrm{Im}\, DF + T_y D = \R^2$ is automatically satisfied, since
$\mathrm{Im}\, DF = \R^2 \supset T_y D$.  Hence
$F \pitchfork D$ for every $D$.

\paragraph{Solution to Exercise~20.8 (Fredholm operators).}
\emph{Right shift:} $L(x_1, x_2, \ldots) = (0, x_1, x_2, \ldots)$.

$\ker L = \{0\}$ (if $L\mathbf{x} = 0$, then $x_n = 0$ for all
$n$), so $\dim \ker L = 0$.

$\mathrm{Im}(L) = \{(y_1, y_2, \ldots) \in \ell^2 : y_1 = 0\}$,
which has codimension~$1$ in $\ell^2$.  Hence
$\mathrm{coker}(L) \cong \R$, with $\dim \mathrm{coker}(L) = 1$.

Since both kernel and cokernel are finite-dimensional, $L$ is
Fredholm, with index $= 0 - 1 = -1$.

\emph{Left shift:} $R(x_1, x_2, \ldots) = (x_2, x_3, \ldots)$.

$\ker R = \{(x_1, 0, 0, \ldots) : x_1 \in \R\} \cong \R$, so
$\dim \ker R = 1$.

$\mathrm{Im}(R) = \ell^2$ (surjective: given any
$(y_1, y_2, \ldots) \in \ell^2$, set $\mathbf{x} = (0, y_1, y_2,
\ldots)$), so $\dim \mathrm{coker}(R) = 0$.

Hence $R$ is Fredholm with index $= 1 - 0 = +1$.

Note that $R = L^*$ (the adjoint of the right shift), and
$\mathrm{index}(L^*) = -\mathrm{index}(L)$, which is consistent:
$+1 = -(-1)$.


\newpage

\begin{thebibliography}{99}

\bibitem{Abraham1963}
R.~Abraham,
Transversality in manifolds of mappings,
\emph{Bull.\ Amer.\ Math.\ Soc.}\ \textbf{69} (1963), 470--474.

\bibitem{Amari1985}
S.-i.~Amari,
\emph{Differential-Geometrical Methods in Statistics},
Lecture Notes in Statistics \textbf{28}, Springer, 1985.

\bibitem{BarndorffNielsen1978}
O.\,E.~Barndorff-Nielsen,
\emph{Information and Exponential Families in Statistical Theory},
Wiley, 1978.

\bibitem{GolubitskyGuillemin1973}
M.~Golubitsky and V.~Guillemin,
\emph{Stable Mappings and Their Singularities},
Graduate Texts in Mathematics \textbf{14}, Springer, 1973.

\bibitem{GuilleminPollack1974}
V.~Guillemin and A.~Pollack,
\emph{Differential Topology},
Prentice-Hall, 1974.

\bibitem{Hirsch1976}
M.\,W.~Hirsch,
\emph{Differential Topology},
Graduate Texts in Mathematics \textbf{33}, Springer, 1976.

\bibitem{JorgensenLabouriau2012}
B.~J{\o}rgensen and R.~Labouriau,
\emph{Exponential Families and Theoretical Inference},
2nd ed., Springer, 2012.

\bibitem{A}
R.~Labouriau,
\emph{Distributional Statistical Models: Weak Moments, Cumulants, and
  a Central Limit Theorem},
arXiv:2604.20634 [math.PR], 2026.

\bibitem{B}
R.~Labouriau,
\emph{Weak Moment Methods for Statistical Inference: with an
  Application to Robust Estimation},
arXiv:2604.23619 [stat.ME], 2026.

\bibitem{C}
R.~Labouriau,
\emph{Statistical Inference Beyond Likelihood via Distributional
  Representations and Estimating Functions},
in preparation, 2026.

\bibitem{QuasiAnalNotes}
R.~Labouriau,
\emph{Notes on Quasi-Analyticity, the Moment Problem, and Weak
  Moments},
lecture notes, 2026.

\bibitem{TransversalityPaper}
R.~Labouriau,
\emph{Transversality and Geometric Regularisation in Distributional
  Statistical Models},
  arXiv:2605.04536 [math.ST], 2026.

\bibitem{Smale1965}
S.~Smale,
An infinite dimensional version of Sard's theorem,
\emph{Amer.\ J.\ Math.}\ \textbf{87} (1965), 861--866.

\bibitem{Stoyanov2000}
J.\,M.~Stoyanov,
Krein condition in probabilistic moment problems,
\emph{Bernoulli} \textbf{6} (2000), 939--949.

\bibitem{Thom1954}
R.~Thom,
Quelques propri\'et\'es globales des vari\'et\'es diff\'erentiables,
\emph{Comment.\ Math.\ Helv.}\ \textbf{28} (1954), 17--86.

\end{thebibliography}
\end{document}